\documentclass[11pt]{amsart}
\usepackage{amssymb,amsmath,amsfonts,amscd,euscript}

\newcommand{\nc}{\newcommand}

\numberwithin{equation}{section}
\newtheorem{thm}{Theorem}[section]
\newtheorem*{thm*}{Theorem}
\newtheorem{prop}[thm]{Proposition}
\newtheorem{lem}[thm]{Lemma}
\newtheorem{exam}[thm]{Example}
\newtheorem{cor}[thm]{Corollary}
\newtheorem*{cor*}{Corollary}
\theoremstyle{remark}
\newtheorem{rem}[thm]{Remark}
\newtheorem{definition}[thm]{Definition}
\newtheorem{dfn}[thm]{Definition}

\usepackage{color}
\definecolor{mycolor}{rgb}{1, 0, 0}
\definecolor{mycolorG}{rgb}{.75,.5,0}
\definecolor{mycolorS}{rgb}{0,0,0}

\newcount\cols {\catcode`,=\active\catcode`|=\active
 \gdef\Young(#1){\hbox{$\vcenter
 {\mathcode`,="8000\mathcode`|="8000
  \def,{\global\advance\cols by 1 &}%
  \def|{\cr
        \multispan{\the\cols}\hrulefill\cr
        &\global\cols=2 }%
  \offinterlineskip\everycr{}\tabskip=0pt
  \dimen0=\ht\strutbox \advance\dimen0 by \dp\strutbox
  \halign
   {\vrule height \ht\strutbox depth \dp\strutbox##
    &&\hbox to \dimen0{\hss$##$\hss}\vrule\cr
    \noalign{\hrule}&\global\cols=2 #1\crcr
    \multispan{\the\cols}\hrulefill\cr%
   }
 }$}}
 \gdef\Skew(#1:#2){\hbox{$\vcenter
 {\mathcode`,="8000\mathcode`|="8000
  \dimen0=\ht\strutbox \advance\dimen0 by \dp\strutbox
  \def\boxbeg{\vbox
    \bgroup\hrule\kern-0.4pt\hbox to\dimen0\bgroup\strut\vrule\hss$}%
  \def\boxend{$\hss\egroup\hrule\egroup}%
  \def,{\boxend\boxbeg}%
  \def|##1:{\boxend\vrule\egroup\nointerlineskip\kern-0.4pt
    \moveright##1\dimen0\hbox\bgroup\boxbeg}%
  \def\\##1\\##2:{\boxend\vrule\egroup\nointerlineskip\kern-0.4pt
    \kern ##1\dimen0\moveright##2\dimen0\hbox\bgroup\boxbeg}%
  \moveright#1\dimen0\hbox\bgroup\boxbeg#2\boxend\vrule\egroup
 }$}} }

\def\smallsquares
 {\textfont0=\scriptfont0 \scriptfont0=\scriptscriptfont0
  \textfont1=\scriptfont1 \scriptfont1=\scriptscriptfont1
  \setbox0=\hbox{$($}
  \setbox\strutbox=\hbox{\vrule width 0pt height\ht0 depth\dp0 }
 }

\nc{\gl}{\mathfrak{gl}}
\nc{\GL}{\mathfrak{GL}}
\nc{\g}{\mathfrak{g}}
\nc{\gh}{\widehat\g}
\nc{\h}{\mathfrak{h}}
\nc{\n}{\mathfrak{n}}
\nc{\la}{\lambda}
\nc{\C}{\mathbb C }
\nc{\D}{\mathbb D }
\nc{\Z}{\mathbb Z }
\nc{\N}{\mathbb N }
\nc{\R}{\mathbb R }
\nc{\Q}{\mathbb Q }
\nc{\al}{\alpha }
\nc{\bs}{{\bf s}}
\nc{\bbs}{\bar{\bf s}}
\nc{\bt}{{\bf t}}
\nc{\br}{{\bf r}}
\nc{\bp}{{\bf p}}
\nc{\bm}{{\bf m}}
\nc{\ba}{\bar{\alpha} }
\nc{\bj}{{\bf j}}
\nc{\bk}{{\bf k}}

\nc{\om}{\omega}
\nc{\ta}{\theta}
\nc{\ve}{\varepsilon}
\nc{\ch}{{\mathop {\rm ch}}}
\nc{\Tr}{{\mathop {\rm Tr}\,}}
\nc{\Id}{{\mathop {\rm Id}}}
\nc{\ad}{{\mathop {\rm ad}}}
\nc{\bra}{\langle}
\nc{\ket}{\rangle}
\nc{\x}{{\bf x}}
\nc{\pa}{\partial}
\nc{\ld}{\ldots}
\nc{\cd}{\cdots}
\nc{\hk}{\hookrightarrow}
\nc{\T}{\otimes}
\newcommand{\bea}{\begin{equation}}
\newcommand{\ena}{\end{equation}}
\newcommand{\be}{\begin{equation*}}
\newcommand{\en}{\end{equation*}}
\nc{\gr}{\mathrm{gr}}
\nc{\ov}{\overline}

\nc{\msl}{\mathfrak{sl}}
\nc{\mgl}{\mathfrak{gl}}
\nc{\U}{\hbox{\rm U}}
\nc{\V}{\EuScript V}

\newcommand{\bc}{{\mathbb C}}
\newcommand{\bz}{{\mathbb Z}}

\newcommand{\fg}{{\mathfrak g}}

\newcommand{\fn}{{\mathfrak n}}

\newcommand{\lam}{\lambda}

\newcommand{\fh}{{\mathfrak h}}
\def\gr{\operatorname{gr}}

\def\wt{\operatorname{wt}}
\def\supp{\operatorname{supp}}

\begin{document}

\title[PBW filtration and bases]
{PBW filtration  and bases for irreducible modules in type ${\tt A}_n$}

\author{Evgeny Feigin, Ghislain Fourier and Peter Littelmann}
\address{Evgeny Feigin:\newline
Tamm Department of Theoretical Physics,
Lebedev Physics Institute,\newline
Leninsky prospect, 53,
119991, Moscow, Russia,\newline
{\it and }\newline
French-Russian Poncelet Laboratory, Independent University of Moscow
}
\email{evgfeig@gmail.com}
\address{Ghislain Fourier:\newline
Mathematisches Institut, Universit\"at zu K\"oln,\newline
Weyertal 86-90, D-50931 K\"oln,Germany
}
\email{gfourier@math.uni-koeln.de}
\address{Peter Littelmann:\newline
Mathematisches Institut, Universit\"at zu K\"oln,\newline
Weyertal 86-90, D-50931 K\"oln,Germany
}
\email{littelma@math.uni-koeln.de}

\begin{abstract}
We study the PBW filtration on the highest weight representations $V(\la)$ of
$\msl_{n+1}$. This filtration is induced by the standard degree filtration on
$U(\n^-)$. We give a description
of the associated graded $S(\n^-)$-module $gr V(\la)$ in terms of generators and
relations. We also construct a basis of $gr V(\la)$. As an application we
derive a graded combinatorial character formula for $V(\la)$, and we obtain
a new class of bases of the modules $V(\la)$ conjectured by Vinberg in 2005. 
\end{abstract}
\maketitle
\section*{Introduction}
Let $\g$ be a simple Lie algebra and let $\g=\n^+\oplus\h\oplus \n^-$ be a Cartan decomposition.
For a dominant integral $\la$ we denote by $V(\la)$ the irreducible $\g$-module with
highest weight $\la$. Fix a highest weight vector $v_\la\in V(\la)$.
Then $V(\la)=\U(\n^-)v_\la$, where $\U(\n^-)$ denotes the universal
enveloping algebra of $\n^-$. The degree filtration $\U(\n^-)_s$ on $\U(\n^-)$ is defined by:
$$
\U(\n^-)_s=\mathrm{span}\{x_1\dots x_l:\ x_i\in\n^-, l\le s\}.
$$
In particular, $\U(\n^-)_0=\C$ and $gr \U(\n^-)\simeq S(\n^-)$, where $S(\n^-)$
denotes the symmetric algebra over $\n^-$. The filtration $\U(\n^-)_s$ induces a filtration
$V(\la)_s$ on $V(\la)$:
$$
V(\la)_s=\U(\n^-)_sv_\la.$$
We call this filtration the
PBW filtration. In this paper we study
the associated graded space $gr V(\la)$ for $\g$ of type ${\tt A}_n$.

So from now on we fix $\g=\msl_{n+1}$. Note that
$gr V(\la)=S(\n^-)v_\la$ is a cyclic $S(\n^-)$-module, so we can write
$$
gr V(\la)\simeq S(\n^-)/I(\la),
$$
for some ideal $I(\la)\subset S(\n^-)$. For example, for any positive root $\al$
the power $f_\al^{(\la,\al)+1}$ of a root vector $f_\al\in\n^-_{-\alpha}$ belongs to $I(\la)$
since $f_\al^{(\la,\al)+1}v_\la=0$ in $V(\la)$. To
describe $I(\la)$ explicitly, consider the action of the opposite subalgebra
$\n^+$ on $V(\la)$. It is easy to see that $\n^+ V(\la)_s\hk V(\la)_s$, so we
obtain the structure of an $U(\n^+)$-module on $gr V(\la)$ as well.
We show:
\vskip 5pt\noindent
{\bf Theorem~A.}\
$I(\la)=S(\n^-)\left(\U(\n^+)\circ \mathrm{span}\{f_\al^{(\la,\al)+1}, \al>0\}\right)$.
\vskip 5pt
Theorem $A$ should be understood as a commutative analogue
of the well-known description of $V(\lam)$ as the quotient
$$
V(\la)\simeq \U(\n^-)/\langle  f_\al^{(\la,\al)+1}, \al>0 \rangle
$$
(see for example \cite{H}).

Our second problem (closely related to the first one)
is to construct a monomial basis of
$gr V(\la)$. The elements $\prod_{\al>0} f_\al^{s_\al} v_\la$ with $s_\al\ge 0$
obviously span $gr V(\la)$ (recall that the order in  $\prod_{\al>0} f_\al^{s_\al}$
is not important since $f_\al$ are considered as elements of $S(\n^-)$).
For each $\la$ we construct a set $S(\la)$ of multi-exponents
$\bs=\{s_\al\}_{\al>0}$ such that the elements
$$
f^\bs v_\la=\prod_{\al>0} f_\al^{s_\al}v_\la, \ \bs\in S(\la)
$$
form a basis of $gr V(\la)$.
To give a precise definition of $S(\la)$ we need to introduce
the notion of a {\it Dyck path}, which occurs already in Vinberg's conjecture: 

Let $\al_1,\dots,\al_n$ be the set of simple roots for $\msl_{n+1}$.
Then all positive roots are of the form $\al_{p,q}=\al_p + \dots + \al_q$
for some $1\le p\le q\le n$.
We call a sequence
\[
\bp=(\beta(0), \beta(1),\dots, \beta(k)), \ k\ge 0,
\]
of positive roots a Dyck path (or simply a path) if
it satisfies the following conditions: either $k=0$, and then $\bp =(\al_i)$ for some simple root
$\alpha_i$, or $k\ge 1$, and then
$\beta(0)=\al_i$, $\beta(k)=\al_j$ for some $1\le i< j\le n$ and
the elements in between obey the following recursion rule:
$$
\text{if } \beta(s)=\al_{p,q} \qquad \text{ then } \qquad \beta(s+1)=\al_{p,q+1} \text{ or }
\beta(s+1)=\al_{p+1,q}.
$$
Denote by $\D$ the set of all Dyck paths.
For a dominant weight $\la=\sum_{i=1}^n m_i\omega_i$
let $P(\lam)\subset \R^{\frac{1}{2}n(n+1)}_{\ge 0}$ be the polytope
\begin{equation}
\label{polytopeequation}
P(\lam):=\bigg\{(r_ \al)_{\al> 0}\mid
\begin{array}{l}
\forall\bp\in\D: \text{ If }\beta(0)=\al_i,\beta(k)=\al_j, \text{ then }\\
r_{\beta(0)} + \dots + r_{\beta(k)}\le m_i + \dots + m_j\\
\end{array}\bigg\},
\end{equation}
and let $S(\la)$ be the set of integral points in $P(\lam)$.
We show:
\vskip 5pt\noindent
{\bf Theorem~B}.\
{\it The set of elements $f^\bs v_\la$, $\bs\in S(\la)$, forms a basis of $gr V(\la)$.}
\vskip 5pt
For $\bs \in S(\lambda)$ define the weight
$$
\wt(\bs) := \sum_{1 \leq j \leq k \leq n} s_{j,k} \alpha_{j,k}.
$$
As an important application we obtain
\begin{cor*}
\begin{itemize}
\item[]
\item[i)] 
For each $\bs\in S(\la)$ fix an arbitrary order of factors $f_\al$ in the product
$\prod_{\al >0} f_\al^{s_\al}$. 
Let $f^\bs=\prod_{\al >0} f_\al^{s_\al}$ be the ordered product.
Then the elements $f^\bs v_\la$, $\bs\in S(\la)$, form a basis of $V(\la)$.
\item[ii)] $\dim V(\lam)=\sharp S(\la)$.
\item[iii)] $char V(\lam) =\sum_{\bs\in S(\lam)} e^{\la-\wt(\bs)}$.
\end{itemize}
\end{cor*}
We note that the order in the corollary above is important since we are back to the
action of the (in general) not commutative enveloping algebra.
We thus obtain a family of bases for irreducible $\msl_{n+1}$-modules.
Motivated by a different background, the existence of these bases (with the same indexing set)
was conjectured by Vinberg (see \cite{V}). Using completely different arguments, he
proved the conjecture for ${\mathfrak{sl}}_4$, for ${\mathfrak{sp}}_4$ and ${\tt G}_2$.
Note also that the data labeling the basis vectors is similar to that
for the Gelfand-Tsetlin patterns (see \cite{GT}). However, these bases are very different
from the GT basis.
\begin{exam}\rm
For $\g=\msl_3(\bc)$, there are only three Dyck paths, the two of length 1 corresponding
to the simple roots, and the path which involves all positive roots. In the following
we write the elements of $P(\lam)$ in a triangular form, where we put $r_1=r_{\alpha_1}$ and
$r_2=r_{\alpha_2}$ in the first row and $r_{12}=r_{\alpha_1+\alpha_2}$ in the second row.
For $\lam=m_1\om_1+m_2\om_2$ the associated polytope is
$$
P(\lam)=\bigg\{\Skew(0:r_1,r_2|.5:r_{12})\in\R_{\ge 0}^3\mid
\begin{array}{l}
0\le r_1\le m_1,0\le r_2\le m_2,\\
r_1+r_2+r_{12}\le m_1+m_2
\end{array}\bigg\},
$$
For the set of integral points we get for example
$$
S(\omega_1)=\bigg\{\Skew(0:0,0|.5: 0),\Skew(0:1,0|.5: 0),\Skew(0:0,0|.5: 1)\bigg\},\
S(\omega_2)=\bigg\{\Skew(0:0,0|.5: 0),\Skew(0:0,1|.5: 0),\Skew(0:0,0|.5: 1)\bigg\}.
$$
and
$$
S(2\omega_1+\omega_2)=\bigg\{\Skew(0:0,0|.5: 0),
\Skew(0:1,0|.5: 0),\Skew(0:2,0|.5: 0),
\Skew(0:0,0|.5: 1),\Skew(0:0,0|.5: 2), \Skew(0:0,0|.5: 3),
\Skew(0:1,0|.5: 1), \Skew(0:1,0|.5: 2),
$$
$$
\Skew(0:2,0|.5: 1),
\Skew(0:0,1|.5: 0),\Skew(0:1,1|.5: 0),\Skew(0:2,1|.5: 0),
\Skew(0:0,1|.5: 1),\Skew(0:0,1|.5: 2),\Skew(0:1,1|.5: 1)\bigg\}.
$$
\end{exam}
\bigskip
We finish the introduction with several remarks. The PBW filtration for representations of
affine Kac-Moody algebras was considered in \cite{FFJMT}, \cite{F1}, \cite{F2}.
It was shown that it has important applications in the representation theory of
current and affine algebras and in mathematical physics.

There exist special representations $V(\la)$ such that the operators
$f^\bs$ consist only of mutually commuting root vectors, even before passing to
$gr V(\la)$. These modules can be described via the theory of abelian radicals
and turned out to be important in the theory of vertex operator algebras
(see \cite{GG}, \cite{FFL}, \cite{FL}).

Let $V_w(\la)\hk V(\la)$ be a Demazure module for some element $w$ from the Weyl group.
For special choices of $w$ there exists a basis of $V_w(\la)$ similar to the one given in Theorem $B$.
We conjecture that this should be true for all $w\in W$ and we will discuss this elsewhere.

Finally we note that $gr V(\la)$ carries an additional grading on each weight space $V(\la)^\mu$
of $V(\la)$:
$$
gr V(\la)^\mu=\bigoplus_{s\ge 0} gr_s V(\la)^\mu=\bigoplus_{s\ge 0} V(\la)^\mu_s/V(\la)^\mu_{s-1}.
$$
The graded character of the weight space is the polynomial
$$
p_{\lam,\mu}(q):=\sum_{s\ge 0}(\dim V(\la)^\mu_s/V(\la)^\mu_{s-1})q^s.
$$
Define the degree $\deg(\bs) := \sum_{1 \leq j \leq k \leq n} s_{j,k}$ for $\bs\in S(\lam)$, and let
$S(\lam)^\mu$ be the subset of elements such that $\mu=\lam-\wt(\bs)$.
Then
\begin{cor*}
$p_{\lam,\mu}(q)=\sum_{\bs\in S(\lam)^\mu} q^{\deg \bs}$.
\end{cor*}
We note that our filtration is different from the Brylinski-Kostant filtration
(see \cite{Br}, \cite{K}).

Our paper is organized as follows:\\
In Section 1 we introduce notations and state the problems.
Sections 2 and 3 are devoted to the proof of Theorem $B$ (see Theorem~\ref{basis}).
In Section 2 we prove the spanning property and in Section 3 the linear independence.
In Section 4 we summarize our constructions and prove Theorem $A$
(see Theorem~\ref{idealtheorem}).

\section{Definitions}
Let $R^+$ be the set of positive roots of $\msl_{n+1}$.  Let
$\al_i$, $\omega_i$ $i=1,\dots,n$ be the simple roots and the fundamental weights.
All roots of $\msl_{n+1}$ are of the form $\al_p + \al_{p+1} +\dots + \al_q$
for some $1\le p\le q\le n$. In what follows we denote such a
root by $\al_{p,q}$, for example $\al_i=\al_{i,i}$.

Let $\msl_{n+1}=\n^+\oplus\h\oplus \n^-$ be the Cartan decomposition.
Consider the increasing degree  filtration on
the universal enveloping algebra of $\U(\n^-)$:
\begin{equation}\label{df}
\U(\n^-)_s=\mathrm{span}\{x_1\dots x_l:\ x_i\in\n^-, l \le s \},
\end{equation}
for example, $\U(\n^-)_0=\C \cdot 1$.

For a dominant integral weight  $\la=m_1\omega_1 + \dots + m_n\omega_n$ let
$V(\la)$ be the corresponding irreducible highest weight
$\msl_{n+1}$-module with a highest weight vector $v_\la$.
Since $V(\la)=\U(\n^-)v_\la$, the filtration \eqref{df} induces an increasing  filtration $V(\la)_s$ on $V(\la)$:
\[
V(\la)_s=\U(\n^-)_s v_\la.
\]
We call this filtration the PBW filtration and study the associated graded
space $gr V(\la)$.  In the following lemma we describe some operators acting on
$gr V(\la)$. Let $S(\n^-)$ denotes the symmetric algebra of $\n^-$.

\begin{lem}
The action of  $\U(\n^-)$ on $V(\la)$ induces the structure of a $S(\n^-)$-module on
$gr V(\la)$ and
\[
gr(V(\la))=S(\n^-)v_\la.
\]
The action of $\U(\n^+)$ on $V(\la)$ induces the structure of a $U(\n^+)$-module on
$gr V(\la)$.
\end{lem}
\begin{proof}
The first statement is obviously true by the definition of the filtrations $\U(\n^-)_s$ and
$V(\la)_s$. The inclusions $\U(\n^+) V(\la)_s\hk V(\la)_s$ imply the second statement.
\end{proof}

Our aims are:
\begin{itemize}
\item to describe $gr V(\la)$ as an  $S(\n^-)$-module, i.e. describe the ideal
$I(\la)\hk S(\n^-)$ such that $gr V(\la)\simeq S(\n^-)/I(\la)$;
\item to find a basis of $gr V(\la)$.
\end{itemize}
The description of the ideal will be given in the last section. To describe the basis
we recall the definition of the Dyck paths:
\begin{dfn}\rm
A {\it Dyck path} (or simply {\it a path}) is a sequence
\[
\bp=(\beta(0), \beta(1),\dots, \beta(k)), \ k\ge 0
\]
of positive roots satisfying the following conditions:
\begin{itemize}
\item[{\it i)}] If $k=0$, then $\bp$ is of the form $\bp =(\al_i)$ for some simple root $\alpha_i$;
\item[{\it ii)}] If $k\ge 1$, then
\begin{itemize}
\item[{\it a)}] the first and last elements are simple roots. More precisely,
$\beta(0)=\al_i$ and $\beta(k)=\al_j$ for some $1\le i< j\le n$;
\item[{\it b)}] the elements in between obey the following recursion rule:
If $\beta(s)=\al_{p,q}$ then the next element in the sequence is of the form either
$\beta(s+1)=\al_{p,q+1}$  or $\beta(s+1)=\al_{p+1,q}.$
\end{itemize}
\end{itemize}
\end{dfn}

\begin{exam}
Here is an example for a path for $\msl_6$:
\[
\bp =(\al_2,\al_2+\al_3,\al_2+\al_3+\al_4,\al_3+\al_4,\al_4,\al_4+\al_5,\al_5).
\]
\end{exam}
\vskip 5pt\noindent
For a multi-exponent $\bs=\{s_ \beta\}_{\beta>0}$, $s_ \beta\in\Z_{\ge 0}$, let $f^\bs$ be the element
\[
f^\bs=\prod_{\beta\in R^+} f_ \beta^{s_ \beta}\in S(\fn^-).
\]
\begin{dfn}
For an integral dominant $\msl_{n+1}$-weight $\la=\sum_{i=1}^n m_i\omega_i$
let $S(\la)$ be the set of all multi-exponents $\bs=(s_ \beta)_{\beta\in R^+}\in\bz_{\ge 0}^{R^+}$
such that for all Dyck paths $\bp =(\beta(0),\dots, \beta(k))$
\begin{equation}\label{upperbound}
s_{\beta(0)}+s_{\beta(1)}+\dots + s_{\beta(k)}\le m_i+m_{i+1}+\dots + m_j,
\end{equation}
where $\beta(0)=\al_i$ and $\beta(k)=\al_j$.
\end{dfn}

In the next two sections we prove the following theorem.

\begin{thm}\label{basis}
The set $f^\bs v_\la$, $\bs\in S(\la)$, form a basis of $gr V(\la)$.
\end{thm}
\proof
In Section 2 we show that the elements $f^\bs v_\la$, $\bs\in S(\la)$, span $gr V(\la)$, see
Theorem~\ref{linearindipencetheorem}. In Section 3 we show that the number $\sharp S(\la)$
is smaller than or equal to $\dim V(\lam)$ (see Theorem~\ref{numberofelements}), which finishes
the proof of Theorem~\ref{basis}.
\qed

\section{The spanning property}\label{span}
The space $gr V(\la)$ is endowed with the structure of a cyclic $S(\n^-)$-module,
i.e. $gr V(\la)=S(\n^-)v_\la$ and hence  $gr V(\la)=S(\n^-)/I(\la)$, where $I(\la)$ is some ideal in $S(\n^-)$.
Our goal in this section is to prove that the elements
$f^\bs v_\la$, $\bs\in S(\la)$, span $gr V(\la)$.

Let $\la=m_1\omega_1 + \dots + m_n\omega_n$. The strategy is as follows:
$f_\al^{(\la,\al)+1}v_\la=0$ in $V(\la)$ for all positive roots $\al$,
so for $\alpha=\al_i + \dots + \al_j$, $i\le j$ we have the relation
\[
f_{\al_i + \dots + \al_j}^{m_i+\dots + m_j+1} \in I(\la).
\]
In addition we have the operators $e_\al$
acting on $gr V(\la)$,
and $I(\la)$ is stable with respect to $e_\al$.
By applying the operators $e_\al$ to $f_{\al_i + \dots + \al_j}^{m_i+\dots + m_j+1}$,
we obtain new relations. We prove that these relations are enough to
rewrite any vector $f^\bt v_\lam$ as a linear combination of $f^\bs v_\lam$ with
$\bs\in S(\la)$.

We start with some notations. For $1\le i < j\le n$ set
\[
\al_{i,j}=\al_i+\dots +\al_j,\ s_{i,j}=s_{\al_{i,j}},\ f_{i,j}=f_{\al_{i,j}},
\]
and for convenience we set $\al_{i,i}=\al_i$, $s_{i,i}=s_{\al_{i}}$ and $f_{i,i}=f_{\al_{i}}$.

By the degree $\deg \bs$ of a multi-exponent we mean the degree of the corresponding
monomial in $S(\fn^-)$, i.e. $\deg \bs=\sum s_{i,j}$.

We are going to define a {\it monomial order} on $S(\fn^-)$. To begin with, we define a total
order on the set of generators $ f_{i,j}$, $1\le i\le j\le n$. We say that $(i,j)\succ (k,l)$ if $i>k$
or if $i=k$ and $j>l$. Correspondingly we say that $f_{i,j}\succ f_{k,l}$ if $(i,j)\succ (k,l)$, so
\[
f_{n,n}\succ f_{n-1,n}\succ f_{n-1,n-1}\succ f_{n-2,n}\succ\ldots \succ f_{2,3}\succ f_{2,2}\succ f_{1,n}\succ\ldots\succ f_{1,1}.
\]
We use the associated {\it homogeneous lexicographic ordering} on the set of monomials in these generators of $S(\fn^-)$.

We use the ``same" total order on the set of multi-exponents, i.e. $\bs\succ \bt$ if and only if $f^\bs\succ f^\bt$.
More explicitly: for two multi-exponents $\bs$ and $\bt$ we write $\bs\succ\bt$:
\begin{itemize}
\item if $\deg \bs>\deg \bt$,
\item if $\deg \bs = \deg \bt$ and there exist $1\le i_0\le j_0\le n$ such
that $s_{i_0j_0}>t_{i_0j_0}$ and for $i>i_0$ and ($i=i_0$ and $j>j_0$) we have $s_{i,j}=t_{i,j}$.
\end{itemize}

\begin{prop}\label{straightening}
Let $\bp=(p(0),\dots,p(k))$ be a Dyck path with $p(0)=\al_i$ and $p(k)=\al_j$.
Let $\bs$ be a multi-exponent supported on $\bp$, i.e. $s_\al=0$ for $\al\notin\bp$.
Assume further that
\[
\sum_{l=0}^k s_{p(l)}>m_i+\dots +m_j.
\]
Then there exist some constants $c_{\bt}$ labeled by multi-exponents $\bt$ such that
\begin{equation}\label{straighteninglaw}
f^\bs +\sum_{\bt<\bs} c_\bt f^\bt \in I(\lam)
\end{equation}
($\bt$ does not have to be supported on $\bp$).
\end{prop}
\begin{rem}\rm
We refer to \eqref{straighteninglaw} as a {\it straightening law} because it implies
$$
f^\bs = - \sum_{\bt<\bs} c_\bt f^\bt \text{\ in\ }S(\fn^-)/I(\lam)\simeq gr V(\lam).
$$
\end{rem}
\begin{proof}
We start with the case $p(0)=\al_1$ and $p(k)=\al_n$ (so, $k=2n-2$).
This assumption is just for convenience. In the general case one
has $\bp$ with $p(0)=\al_i$, $p(k)=\al_j$ and one would start with the relation
$f_{i,j}^{m_i+\dots +m_j+1} \in I(\lam)$ instead of the relation
$f_{1,n}^{m_1+\dots +m_n+1} \in I(\lam)$ below.

So from now on we assume without loss of generality that $p(0)=\al_1$ and $p(k)=\al_n$.
Let $S_+(\fh\oplus \fn^+)\subset S(\fh\oplus \fn^+)$ be the maximal homogeneous ideal of polynomials
without constant term. The adjoint action of $U(\fn^+)$ on
$\fg$ induces an action of $U(\fn^+)$ on $S(\fg)$ and hence on
$$
S(\fn^-)\simeq S(\fg)/S(\fn^-)S_+(\fh\oplus\fn^+).
$$
In the following we use the differential operators $\pa_\al$ defined by
$$
\pa_\al f_\beta=
\begin{cases}
f_{\beta-\al},\  \text{ if }  \beta-\al\in\triangle^+,\\
0,\ \text{ otherwise}.
\end{cases}
$$
The operators $\pa_\al$ satisfy the property
\[
\pa_\al f_\beta = c_{\al,\beta} (\text{ad}\, e_\al)(f_\beta),
\]
where $c_{\al,\beta}$ are some non-zero constants.
In the following we use very often the following consequence:
if  $f_{\beta_1}\dots f_{\beta_l}\in I(\la)$, then for any $\al_1,\dots,\al_s$
\[
\pa_{\al_1}\dots\pa_{\al_s} f_{\beta_1}\dots f_{\beta_l}\in I(\la).
\]
Since $f_{1,n}^{m_1+\dots +m_n+1} v_\lam=0$
in $gr V(\la)$, it follows that
\[
f_{1,n}^{s_{p(0)}+\dots + s_{p(k)}} \in I(\lam).
\]
Write $\pa_{i,j}$ for $\pa_{\al_{i,j}}$. For instance, we have
\begin{equation}\label{pa1}
\pa_{1,i} f_{1,j}=f_{i+1,j},\ \pa_{j,n}f_{i,n}=f_{i,j-1} \ \text{ for } 1\le i<j\le n.
\end{equation}
For $i,j=1,\dots,n$ set
\[
s_{\bullet, j}=\sum_{i=1}^j s_{i,j},\quad s_{i,\bullet}=\sum_{j=i}^n s_{i,j}.
\]
We consider first the vector
\begin{equation}\label{pape1}
\pa_{n,n}^{s_{\bullet, n-1}}\pa_{n-1,n}^{s_{\bullet, n-2}}\dots \pa_{2,n}^{s_{\bullet, 1}}
f_{1,n}^{s_{p(0)}+\dots + s_{p(k)}} \in I(\lam).
\end{equation}
Because of the formulas in \eqref{pa1} we get:
\[
\pa_{2n}^{s_{\bullet, 1}}
f_{1,n}^{s_{p(0)}+\dots + s_{p(k)}}=c_1 f_{1,n}^{s_{p(0)}+\dots + s_{p(k)-s_{\bullet, 1}}}f_{1,1}^{s_{\bullet, 1}}
\]
for some nonzero constant $c_1$, and
\[
\pa_{3n}^{s_{\bullet, 2}}\pa_{2n}^{s_{\bullet, 1}}
f_{1,n}^{s_{p(0)}+\dots + s_{p(k)}} =c_2 f_{1,n}^{s_{p(0)}+\dots + s_{p(k)}-s_{\bullet 1}-s_{\bullet 2}}f_{1,1}^{s_{\bullet 1}}f_{1,2}^{s_{\bullet 2}}
\]
for some nonzero constant $c_2$ etc.
Summarizing, the vector \eqref{pape1} is proportional (with a nonzero constant) to
\[
f_{1,1}^{s_{\bullet, 1}}f_{1,2}^{s_{\bullet ,2}}\dots f_{1,n}^{s_{\bullet, n}} \in I(\lam).
\]
To prove the proposition, we apply more differential operators to the monomial
$f_{1,1}^{s_{\bullet ,1}}f_{1,2}^{s_{\bullet, 2}}\dots f_{1,n}^{s_{\bullet, n}}$. Consider the following element in $I(\lam)\subset S(\fn^-)$:
\begin{equation}\label{Aoperator}
A=\pa_{1,1}^{s_{2,\bullet}}\pa_{1,2}^{s_{3,\bullet}}\dots \pa_{1,n-1}^{s_{n,\bullet}}
f_{1,1}^{s_{\bullet ,1}}f_{1,2}^{s_{\bullet, 2}}\dots f_{1,n}^{s_{\bullet,n}}.
\end{equation}
{\bf We claim:}
\begin{equation}\label{Aoperatorequation}
A =\sum_{\bt\le \bs} c_\bt f^\bt \text{\ for some $c_\bs\ne 0$}.
\end{equation}
Now $A \in I(\lam)$ by construction, so the claim proves the proposition.
\vskip 5pt\noindent
{\it Proof of the claim:}
In order to prove the claim we need to introduce some more notation. For $j=1,\ldots,n-1$ set
\begin{equation}\label{Ajoperator}
A_j=\pa_{1,j}^{s_{j+1,\bullet}}\pa_{1,j+1}^{s_{j+2,\bullet}}\dots \pa_{1,n-1}^{s_{n,\bullet}}
f_{1,1}^{s_{\bullet, 1}}f_{1,2}^{s_{\bullet, 2}}\dots f_{1,n}^{s_{\bullet, n}},
\end{equation}
so $A_1=A$. To start an inductive procedure, we begin with $A_{n-1}$:
$$
A_{n-1}=\pa_{1,n-1}^{s_{n,\bullet}} f_{1,1}^{s_{\bullet, 1}}f_{1,2}^{s_{\bullet, 2}}\dots f_{1,n}^{s_{\bullet, n}}.
$$
Now $s_{n,\bullet}=s_{n,n}$ and $\pa_{1,n-1}f_{1,j}=0$ except for $j=n$, so
\begin{equation}\label{Aoperatorequation1}
A_{n-1}= c
f_{1,1}^{s_{\bullet, 1}}f_{1,2}^{s_{\bullet, 2}}\dots f_{1,n}^{s_{\bullet, n}-s_{n,n}}f_{n,n}^{s_{n,n}},
\end{equation}
for some nonzero constant $c$.

The proof will now proceed by decreasing induction. Since the induction procedure
is quite involved and the initial step does not reflect the problems
occurring in the procedure, we discuss for convenience the case $A_{n-2}$
separately.

Consider $A_{n-2}$, up to a nonzero constant we have:
$$
A_{n-2}=\pa_{1,n-2}^{s_{n-1,\bullet}}
f_{1,1}^{s_{\bullet, 1}}f_{1,2}^{s_{\bullet, 2}}\dots f_{1,n}^{s_{\bullet, n}-s_{n,n}}f_{n,n}^{s_{n,n}}.
$$
Now $\pa_{1,n-2}f_{1,j}=0$ except for $j=n-1,n$, and $\pa_{1,n-2}f_{n,n}=0$, so
$$
A_{n-2}=\sum_{\ell=0}^{s_{n-1 ,\bullet}} c_\ell
f_{1,1}^{s_{\bullet, 1}}f_{1,2}^{s_{\bullet, 2}}\dots f_{1,n-1}^{s_{\bullet, n-1}-s_{n-1,\bullet}+\ell}f_{1,n}^{s_{\bullet, n}-s_{n,n}-\ell}
f_{n-1,n-1}^{s_{n-1,\bullet}-\ell}f_{n-1,n}^{\ell}f_{n,n}^{s_{n,n}}.
$$
We need to control which powers $f_{n-1,n}^{\ell}$ can occur. Recall that $\bs$ has support in $\bp$.
If $\al_{n-1}\not\in \bp$, then $s_{n-1,n-1}=0$ and
$s_{n-1,\bullet}=s_{n-1,n}$, so $f_{n-1,n}^{s_{n-1,n}}$ is the highest power occurring in the sum.
Next suppose $\al_{n-1}\in \bp$. This implies $\al_{j,n}\not \in \bp$ unless $j=n-1$ or $n$. Since
$\bs$ has support in $\bp$, this implies
$$
s_{\bullet, n}=s_{1,n}+\ldots + s_{n-1,n}+s_{n,n}=s_{n-1,n}+s_{n,n},
$$
and hence again the highest power of $f_{n-1,n}$ which can occur is $f_{n-1,n}^{s_{n-1,n}}$,
and the coefficient is nonzero. So we can write
\begin{equation}\label{Aoperatorequation2}
A_{n-2}=\sum_{\ell=0}^{s_{n-1,n}} c_\ell
f_{1,1}^{s_{\bullet, 1}}
\dots f_{1,n-1}^{s_{\bullet, n-1}-s_{n-1,\bullet}+\ell}f_{1,n}^{s_{\bullet, n}-s_{n,n}-\ell}
f_{n-1,n-1}^{s_{n-1,\bullet}-\ell}f_{n-1,n}^{\ell}f_{n,n}^{s_{n,n}}.
\end{equation}
For the inductive procedure we make the following assumption:

$A_{j}$ is a sum of monomials of the form
\begin{equation}\label{Aoperatorequation3}
\underbrace{f_{1,1}^{s_{\bullet, 1}}\ldots f_{1,j}^{s_{\bullet, j}}f_{1,j+1}^{s_{\bullet, j+1}-*}\ldots f_{1,n}^{s_{\bullet, n}-*}}_X
\underbrace{f_{j+1,j+1}^{t_{j+1,j+1}}f_{j+1,j+2}^{t_{j+1,j+2}}\dots f_{n-1,n}^{t_{n-1,n}}f_{n,n}^{t_{n,n}}}_Y
\end{equation}
having the following properties:
\begin{itemize}
\item[{\it i)}] With respect to the homogeneous lexicographic ordering, all the multi-exponents of the summands, except one,
are strictly smaller than $\bs$.
\item[{\it ii)}] More precisely, there exists a pair $(k_0,\ell_0)$ such that $k_0\ge j+1$,
$s_{k_0 \ell_0}>t_{k_0 \ell_0}$ and $s_{k\ell}=t_{k\ell}$ for all $k>k_0$ and all pairs $(k_0.\ell)$ such that $\ell>\ell_0$.
\item[{\it iii)}] The only exception is the summand such that $t_{\ell,m}=s_{\ell,m}$ for all $\ell\ge j+1$
and all $m$.
\end{itemize}
The calculations above show that this assumption holds for $A_{n-1}$ and $A_{n-2}$.

We come now to the induction procedure and we consider $A_{j-1}= \pa_{1,j-1}^{s_{j,\bullet}}A_{j}$.
Note that $\pa_{1,j-1} f_{1,\ell}=0$ except for $\ell \ge j$, and in this case we have $\pa_{1,j-1} f_{1,\ell}=f_{j,\ell}$.
Furthermore, $\pa_{1,j-1} f_{k,\ell}=0$ for $k\ge j+1$, so applying $\pa_{1,j-1}$ to a summand of the form in (\ref{Aoperatorequation3})
does not change the $Y$-part in (\ref{Aoperatorequation3}). Summarizing,
applying $\pa_{1,j-1}^{s_{j,\bullet}}$ to a summand of the form in (\ref{Aoperatorequation3}) gives a sum of
monomials of the form
\begin{equation}\label{Aoperatorequation4}
\underbrace{f_{1,1}^{s_{\bullet ,1}}\ldots f_{1,j-1}^{s_{\bullet, j-1}}f_{1,j}^{s_{\bullet,j}-*}\ldots f_{1,n}^{s_{\bullet,n}-*}}_{X'}
\underbrace{f_{j,j}^{t_{j,j}}\ldots f_{j,n}^{t_{j,n}} }_Z
\underbrace{f_{j+1,j+1}^{t_{j+1,j+1}}f_{j+1,j+2}^{t_{j+1,j+2}}\dots
f_{n,n}^{t_{n,n}}}_Y.
\end{equation}
We have to show that these summands satisfy again the conditions ${\it i)}$--${\it iii)}$ above (but now for the index $(j-1)$). If we start
in (\ref{Aoperatorequation3}) with a summand which is not the maximal summand, but such that {\it i)} and {\it ii)} hold for the index $j$,
then the same holds obviously also for the index $(j-1)$ for all summands in (\ref{Aoperatorequation4}) because the $Y$-part
remains unchanged.

So it remains to investigate the summands of the form (\ref{Aoperatorequation4}) obtained by applying $\pa_{1j-1}^{s_{j,\bullet}}$ to the only
summand in (\ref{Aoperatorequation3}) satisfying {\it iii)}.

To formalize the arguments used in the calculation for $A_{n-2}$ we need the following notation.
Let $1\le k_1\le k_2\le\dots\le k_n\le n$ be numbers defined by
\[
k_i=\max\{j:\ \al_{i,j}\in\bp\}.
\]
For convenience we set $k_0=1$.
\begin{exam}
For $\bp=(\al_{11},\al_{12},\dots,\al_{1n},\al_{2n},\dots,\al_{n,n})$ we have
$k_i=n$ for all  $i=1,\ldots,n$.
\end{exam}
%
Since $\bs$ is supported on $\bp$ we have
\begin{equation}\label{sumpath}
s_{i,\bullet}=\sum_{\ell=k_{i-1}}^{k_i} s_{i,\ell},\ s_{\bullet, \ell}=\sum_{i:\ k_{i-1}\le \ell\le k_i} s_{i,\ell}.
\end{equation}
Suppose now that we have a summand of the form in (\ref{Aoperatorequation4}) obtained by applying
$\pa_{1j-1}^{s_{j,\bullet}}$ to the only summand in (\ref{Aoperatorequation3}) satisfying {\it iii)}.
Since the $Y$-part remains unchanged, this implies already $t_{n,n}=s_{n,n},\ldots, t_{j+1,j+1}=s_{j+1,j+1}$.
Assume that we have already shown $t_{j,n}=s_{j,n},\ldots, t_{j,\ell_0+1}=s_{j,\ell_0+1}$, then we have to show
that  $t_{j,\ell_0}\le s_{j,\ell_0}$.

We consider five cases:
\begin{itemize}
\item $\ell_0 > k_{j}$. In this case the root $\al_{j,\ell_0}$ is not in the support of $\bp$ and hence
$s_{j,\ell_0}=0$. Since $\ell_0>k_{j}\ge k_{j-1}\ge\ldots\ge k_1$, for the same reason we have
$s_{i,\ell_0}=0$ for $i\le j$.
Recall that the power of $f_{1,\ell_0}$ in $A_{j-1}$ in (\ref{Ajoperator}) is equal to $s_{\bullet, \ell_0}$.
Now $s_{\bullet, \ell_0}=\sum_{i>j} s_{i,\ell_0}$ by the discussion above, and hence $f_{1,\ell_0}^{s_{\bullet, \ell_0}}$
has already been transformed completely by the operators $\pa_{1,i}$, $i>j$,
and hence $t_{j,\ell_0}=0=s_{j,\ell_0}$.
\item $k_{j-1}<\ell_0\le k_{j}$. Since $\ell_0 > k_{j-1}\ge\ldots\ge k_1$, for the same reason as above
we have $s_{i,\ell_0}=0$ for $i<j$, so $s_{\bullet, \ell_0}=\sum_{i\ge j} s_{i,\ell_0}$.
The same arguments as above show that  for the operator $\pa_{1,j-1}$
only the power $f_{1,\ell_0}^{s_{j,\ell_0}}$ is left to be transformed into a power of $f_{j,\ell_0}$,
so necessarily  $t_{j,\ell_0}\le s_{j,\ell_0}$.
\item $k_{j-1}=\ell_0=k_{j}$. In this case $s_{j,\bullet}=s_{j,\ell_0}$ and thus
the operator $\pa_{1,j-1}^{s_{j,\bullet}}=\pa_{1,j-1}^{s_{j,\ell_0}}$ can transform
a power $f^*_{1,\ell_0}$ in $A_j$ only into a power $f^q_{j,\ell_0}$ with $q$ at
most $s_{j,\ell_0}$.
\item $k_{j-1}= \ell_0 < k_{j}$. In this case $s_{j,\bullet}=s_{j ,\ell_0}+s_{j ,\ell_0+1}+\ldots+s_{j ,k_j}$.
Applying $\pa_{1,j-1}^{s_{j,\bullet}}$ to the only summand in (\ref{Aoperatorequation3}) satisfying {\it iii)},
the assumption $t_{j,n}=s_{j,n},\ldots, t_{j,\ell_0+1}=s_{j,\ell_0+1}$ implies that one has to
apply $\pa_{1,j-1}^{s_{j,k_j}}$ to $f_{1,k_j}^*$ and $\pa_{1,j-1}^{s_{j,k_j-1}}$  to $f_{1,k_j-1}^*$ etc. to get the demanded
powers of the root vectors. So for $f^*_{1,\ell_0}$ only the operator $\pa_{1,j-1}^{s_{j,\ell_0}}$ is left for transformations
into a power of $f_{j,\ell_0}$
and hence $t_{j,\ell_0}\le s_{j,\ell_0}$.
\item $\ell_0< k_{j-1}$. In this case $s_{j,\ell_0}=0$ because the root is not in the support.
Since $t_{j,\ell }=s_{j,\ell}$ for $\ell>\ell_0$ and $s_{j,\ell}=0$ for $\ell \le \ell_0$ (same reason as above)
we obtain
\[
\pa_{1,j-1}^{s_{j,\bullet}}=\pa_{1,j-1}^{\sum_{\ell>\ell_0} s_{j,\ell}}.
\]
But by assumption we know that $\pa_{1,j-1}^{s_{j,\ell}}$ is needed to transform
the power  $f_{1,\ell}^{s_{j,\ell}}$ into  $f_{j,\ell}^{s_{j,\ell}}$ for all $\ell >\ell_0$,
so no power of $\pa_{1,j-1}$ is left and thus $t_{j,\ell_0}=0=s_{j,\ell_0}$.
\end{itemize}
It follows that all summands except one satisfy the conditions {\it i),ii)} above.
The only exception is the term where the powers of the operator $\pa_{1,j-1}^{s_{j,\bullet}}$
are distributed as follows:
$$
f_{1,1}^{s_{\bullet, 1}}... f_{1,j-1}^{s_{\bullet, j-1}}
(\pa_{1,j-1}^{s_{j,j}}f_{1,j}^{s_{\bullet, j}}) (\pa_{1,j-1}^{s_{j,j+1}}f_{1,j+1}^{s_{\bullet, j+1}-*}) ...
(\pa_{1,j-1}^{s_{j,n}} f_{1,n}^{s_{\bullet, n}-*}) f_{j+1,j+1}^{s_{j+1,j+1}}... f_{n,n}^{s_{n,n}}.
$$
By construction, this term is nonzero and satisfies the condition {\it iii)}, which finishes
the proof of the proposition.
\end{proof}

\begin{thm}\label{linearindipencetheorem}
The elements $f^{\bs} v_\lambda$ with $\bs\in S(\la)$ span the module $gr V(\la)$.
\end{thm}
\begin{proof}
The elements $f^\bs v_\la$, $\bs$ arbitrary multi-exponent, span $S(\fn^-)/I(\lam)\simeq gr V(\la)$.
We use now the equation \eqref{straighteninglaw} in Proposition~\ref{straightening} as a
straightening algorithm to express $f^\bt v_\la $, $\bt$ arbitrary,
as a linear combination of elements $f^\bs v_\la$ such that $\bs\in S(\lam)$.

Let $\lam=\sum_{i=1}^n m_i\omega_i$ and suppose $\bs\notin S(\la)$, then there exists a Dyck path $\bp=(p(0),\dots,p(k))$
with $p(0)=\al_i$, $p(k)=\al_j$ such that
\[
\sum_{l=0}^k s_{p(l)} > m_i +\dots + m_j.
\]
We define a new multi-exponent $\bs'$ by setting
\[\bs'_\al=
\begin{cases}
s_\al, \ \al\in\bp,\\
0,\ otherwise.
\end{cases}
\]
For the new multi-exponent $\bs'$ we still have
$$
\sum_{l=0}^k s'_{p(l)} > m_i +\dots + m_j.
$$
We can now apply Proposition~\ref{straightening} to $\bs'$ and conclude
$$
f^{\bs'} =\sum_{\bs'>\bt'} c_{\bt'} f^{\bt'}\quad\text{in}\quad S(\fn^-)/I(\lam).
$$
We get $f^\bs$ back as $f^\bs=f^{\bs'}\prod_{\beta\notin\bp} f_\beta^{s_\beta}$.
For a multi-exponent $\bt'$ occurring in the sum with $c_{\bt'}\not=0$ set
$f^\bt= f^{\bt'}\prod_{\beta\notin\bp} f_\beta^{s_\beta}$ and $c_\bt=c_{\bt'}$. Since we have a monomial
order it follows:
\begin{equation}\label{tprime3}
f^\bs =f^{\bs'}\prod_{\beta\notin\bp} f_\beta^{s_\beta}=\sum_{\bs>\bt} c_\bt f^\bt \quad\text{in}\quad S(\fn^-)/I(\lam).
\end{equation}
The equation \eqref{tprime3} provides an algorithm
to express $f^\bs$ in $S(\fn^-)/I(\lam)$ as a sum of elements of the desired form: if some
of the $\bt$ are not elements of $S(\lam)$, then we
can repeat the procedure and express the $f^\bt$ in $S(\fn^-)/I(\lam)$ as a sum of
$f^\br$ with $\br<\bt$. For the chosen ordering any strictly decreasing
sequence of multi-exponents is finite, so after a finite number of steps one obtains an
expression of the form $f^\bs=\sum c_\br f^\br$ in $S(\fn^-)/I(\lam)$ such that
$\br\in S(\lam)$ for all $\br$.
\end{proof}

\section{The linear independence}\label{li}
In the following let $R_i$ denote the subset
$$
R_i = \{ \alpha \in R^+ \; | \; (\omega_i, \alpha) = 1 \}.
$$
We define for a dominant weight $\lambda \in P^+$
$$
R_{\lambda} = \{ \alpha \in R^+ \; | \; (\lambda, \alpha) > 0\}.
$$
Recall that we use $\al_{i,j}$ as an abbreviation for $\al_i+\al_{i+1}+\ldots+\al_j$
(see Section \ref{span}). The set $R_i$ can then be described as
$$
R_i = \{ \alpha_{j,k} \; | \; 1 \leq j \leq i \leq k \leq n \}.
$$
We say a path $\bp$ has {\it color} $i$ if $\exists \; j$ s.t. $\beta(j) \in R_i$. Note that a path
can have several different {colors}.

To simplify the notation we often just write $(j,k)$ for the root $\alpha_{j,k}$ (if no confusion is possible).

Let $\lambda = \sum_{j=1}^{n} m_j \omega_j$ and let $i$ be minimal such that $m_i \neq 0$.
For $\bs \in S(\lambda)$, we denote
$$R^{\bs}_i = \{ (j,k) \in R_i \; | \; s_{j,k} \neq 0 \}.$$
We define two different orders on $R$, a partial order ``$\leq$":
$$
(j_1, k_1) \leq (j_2, k_2) \Leftrightarrow (j_1 \leq j_2 \; \wedge k_1 \leq k_2),
$$
and a total order ``$\ll$'':
$$(j_1, k_1) \ll (j_2, k_2) \Leftrightarrow \text{ if }(k_1 < k_2) \text{ or } (k_1 = k_2 \wedge j_1 < j_2).$$
By definition, ``$\ll$'' covers ``$\leq$''.\\

\begin{exam}\rm
For $\g=\msl_4$ and $i =2$, the minimal element of $R_i$ with respect to both orders is
$(1,2) = \alpha_{1,2} = \alpha_1 + \alpha_2$. Note that $ \alpha_1 + \alpha_2 + \alpha_3 \ll \alpha_2$,
but the two are not comparable with respect to ``$\leq$".
\end{exam}
A tuple $\bs \in S(\lambda)$ will be considered as an ordered tuple with respect to the order ``$\ll$":
$$
\bs = (s_{1,1}, s_{1,2}, s_{2,2}, s_{1,3}, s_{2,3}, s_{3,3}, \ldots , s_{n, n}).
$$
The induced lexicographic order on $S(\lambda)$ is a total order which we again
denote by ``$\ll$".
\begin{rem}\rm
The total order $\ll$ is different from the order $\prec$ used in Section~\ref{span}.
\end{rem}
\begin{exam}\rm
For $\g=\msl_4$ let $\bs$ be defined by
$$
s_{13} = 1, s_{22} = 1 \text{ and } s_{j,k} = 0 \text{ otherwise,}
$$
and let $\bt$ be defined by
$$
t_{12} = 1, t_{23} = 1 \text{ and }  t_{j,k} = 0 \text{ otherwise. }
$$
Then $\bs= (0,0,1,1,0,0)$ and $\bt=(0,1,0,0,1,0)$, and so $\bs \ll \bt$.
\end{exam}
\begin{definition} For $\bs \in S(\lambda)$
denote by $M^{\bs}_i$ the {\it set of minimal elements} in $R^{\bs}_i$ with respect to $\leq$.
We denote by $\bm^{\bs}_i$ the tuple $m_{j,k} = 1$ if $(j,k) \in M^{\bs}_i$ and $m_{j,k} = 0$ otherwise.
\end{definition}
\begin{exam}\rm
1) If $R^{\bs}_i=R_i$, then $M^{\bs}_i=\{\al_{1,i}\}$.
\par
2) If $R^{\bs}_i=\{\al_{i,i},\al_{i-1,i+1},\ldots,\al_{i-\ell,i+\ell}\}$ for some $\ell\le i$, then
$M^{\bs}_i=R^{\bs}_i$.
\end{exam}

\begin{rem}\label{minimalelements}
$1)$.\ For any multi-exponent $\bs$ we have
$$
M^{\bs}_i=\{ \al_{j_l, k_l} \; | \; l=1, \ldots, m \}
$$
for some $m$, and the indices have the property
$$
1\le j_1 < j_2 < \ldots < j_m\le i \le k_m< \dots <  k_2 < k_1 \le n.
$$
If $\bs\in S(\omega_i)$, then
for the associated tuple $\bm_i^\bs$ we get: $\bm_i^\bs=\bs$.\\
\noindent $2)$.\ The sets $M^{\bs}_i$ satisfy the following important property:
any Dyck path contains at most one element of $M^{\bs}_i$, because
the elements of a Dyck path are linearly ordered with respect to ``$\ge$''.
\end{rem}

\begin{prop}\label{minimalset}
For $\bs \in S(\lambda)$ let $M^{\bs}_i$ be the minimal set.
Then $\bm^{\bs}_i \in S(\omega_i)$, and
if $\bs'$ is such that $\bs = \bs' + \bm^{\bs}_i$, then $\bs' \in S(\lambda - \omega_i)$.
\end{prop}
\proof
Note that $\bm^{\bs}_i \in S(\omega_i)$ by Remark~\ref{minimalelements}.
Let $\bs'$ be such that $\bs = \bs' + \bm^{\bs}_i$.
We claim that $\bs' \in S(\lambda - \omega_i)$.
Let $\lam=\sum_{j=i}^n m_j\om_j$. For a Dyck path $\bp$ let
$q_\bp^\lam=\sum_{j\, \text{color of}\, \bp} m_j$
be the upper bound for the defining inequality \eqref{upperbound}
of $S(\lam)$ associated to $\bp$.

If $\bp$ is a Dyck path such that $i$ is not a color, then $q_\bp^\lam= q_\bp^{\lam-\om_i}$
and $s_{\beta}=s_{\beta}'$ for $\beta\not\in R_i$, so $\bs'$ satisfies the defining
inequality for $S(\lambda - \omega_i)$ given by $\bp$.

Let $\bp$ be a Dyck path of color $i$, so $q_\bp^{\lam-\om_i}=q_\bp^\lam-1$. If
$\bp\cap M^{\bs}_i\not=\emptyset$, then $\sum_{(j,k) \in \bp} s_{j,k}' =
\sum_{(j,k) \in \bp} s_{j,k}-1\le q_\bp^\lam-1=q_\bp^{\lam-\om_i}$, so $\bs'$
satisfies the defining inequality for $S(\lambda - \omega_i)$ given by $\bp$.

Suppose now that $\bp$ is a Dyck path of color $i$ but $\bp\cap M^{\bs}_i=\emptyset$.
Recall that the elements in $\supp \bp$ are linearly ordered. Let
$\al_{l,m}$ be the minimal element in $R^{\bs}_i \cap \supp \bp$.
Since $i$ is minimal such that $m_i>0$, note that $s_\beta=0$ for all
$\beta\in\supp\bp$ be such that $\beta<\al_{l,m}$. By assumption,
$\al_{l,m}\not\in M^{\bs}_i$, so let $\al_{r,t}\in M^{\bs}_i$
such that $\al_{r,t}<\al_{l,m}$. Let $\tilde\bp$ be the Dyck path
$$
(\al_{r,r},\al_{r,r+1},\ldots,\al_{r,t},\al_{r,t+1},\ldots,\al_{r,m},\al_{r+1,m},\ldots,\al_{l,m},\beta_1,\ldots,\beta_N),
$$
where $\{\beta_1,\ldots,\beta_N\}$ are the elements in $\supp\bp$ such that $\beta_j> \al_{l,m}$.
Since $\al_{r,t}\in \supp\tilde\bp$ we know:
$$
\sum_{(j,k) \in \bp} s_{j,k}<\sum_{(j,k) \in \tilde \bp} s_{j,k}\le q_\bp^\lam
$$
and hence $\sum_{(j,k) \in \bp} s_{j,k}=\sum_{(j,k) \in \bp} s'_{j,k}\le q_\bp^\lam-1=q_\bp^{\lam-\om_i}$.
\qed\\

For $\bs \in S(\lambda)$ we define a \textit{mutation} of $\bs$ as follows:
\begin{definition}
Let
$$\beta = \sum_{(j,k) \in R} s_{j,k} \alpha_{j,k}$$
and suppose
$$\beta = \sum_{(j,k) \in R} t_{j,k} \alpha_{j,k}$$
where
$$t_{j,k} = 0 \text{ if } (j,k) \notin R_{\lambda} \; ; \; t_{j,k} \geq 0 \text{ if } (j,k) \in R_{\lambda},$$
for some $\bt = (t_{j,k}) \notin S(\lambda)$. Then we call $\bt$ a mutation of $\bs$.
\end{definition}
\begin{exam}\rm
Let $\fg = \mathfrak{sl}_3$ and $\lambda = \omega_2$. Define
$$\bs \text{ by } s_{1,3} = 1, s_{2,2} = 1 \text{ and } s_{i,j} = 0 \text{ else, }$$
and $$\bt \text{ by } t_{1,2} = 1, t_{2,3} = 1 \text{ and } t_{i,j} = 0 \text{ else.}$$
Then $\bt$ is a mutation of $\bs$.
\end{exam}
\begin{prop}\label{mutation}
For $\bs \in S(\lambda)$ let $M^{\bs}_i$ be the minimal set.
If $\bt^1$ is a mutation of $\bm^{\bs}_i$, $\bt = \bt^2 + \bt^1 \in S(\lambda)$
and $t^{1}_{j,k} \geq 0$, then $\bm^{\bs}_i \ll\bm^{\bt}_i$.
\end{prop}
\proof
Recall (see Remark~\ref{minimalelements}) that
$M^{\bs}_i=\{ (j_l, k_l) \; | \; l=1, \ldots, m \}$ with
$$
1\le j_1 < \dots < j_m\le i\le k_m < \dots < k_1.
$$
Let $\bt^1$ be a mutation of $\bm^{\bs}_i$, so $t^1_{j,k}= 0$ for $(j,k) \notin R_i$.
Then there exists $\sigma \in S_{m} \setminus\{ id \}$ such that if
$t^1_{p,q} \neq 0$, then $(p,q) = (j_l, k_{\sigma(l)})$ for
some $1 \leq l \leq m$.
We can even assume that $\sigma(l) \neq l$ for all $l$, because
otherwise $(j_l, k_l)$ is not mutated and appears in
$\bm^{\bs}_i$ and $\bt^1$.\\
It is clear that $\bm^{\bt^1}_i \ll \bm^{\bt}_i$ (or equal), so it suffices to show that $\bm^{\bs}_i \ll \bm^{\bt_1}_i$.\\
Let $x = \sigma^{-1}(m)$, we claim that $M^{\bt}_i \subset \{ (j_1, k_{\sigma(1)}), \ldots, (j_x, k_{\sigma(x)}) \}$.
Let $l > x$, then $j_x < j_l$ and $k_m > k_{\sigma(l)}$ (since $\sigma(l) \neq m$). So $(j_x, k_m) < (j_l, k_{\sigma(l)})$
for all $l > x$.
\qed

\begin{thm}\label{numberofelements}
Let $\lambda = \sum_j m_j \omega_j \in P^+$.
For each $\bs\in S(\la)$ fix an arbitrary order of factors $f_\al$ in the product
$\prod_{\al>0} f_\al^{s_\al}$.
Let $f^\bs=\prod_{\al>0} f_\al^{s_\al}$ be the ordered product in $U(\n^-)$.
Then the elements $f^\bs v_\la$, $\bs\in S(\la)$, form a basis of $V(\la)$.
\end{thm}
\proof
We will prove the claim by induction on $m = \sum_{j=1}^n m_j$.
By Theorem~\ref{linearindipencetheorem} we know that the $f^\bs v_{\lam}$
span the representation $V(\lam)$, so $\dim V(\lam)\ge \sharp S(\lam)$.
For the initial step $m=1$ the description of $S(\om_i)$ in Remark~\ref{minimalelements}
shows that the tuples have all different weights and hence the $f^\bs v_{\om_i}$
are also linearly independent, which proves the claim for the
fundamental representations.

We assume that the claim holds for $\lam$, we want to prove it for $\lambda + \omega_i$.
We may assume again that $i$ is minimal such that $m_i \neq 0$.
The highest weight vector $v_{\lambda} \otimes v_{\omega_i}$ generates
$V(\lambda + \omega_i) \subset V(\lambda) \otimes V(\omega_i)$.
We assume in the following that the roots are ordered in such a way that the
$f_{\alpha}$ with $\alpha \in R_i$ are at the beginning.
Every element $\bs \in S(\lambda + \omega_i)$
defines a vector of $f^{\bs} (v_{\lambda} \otimes v_{\omega_i}) \in V(\lambda + \omega_i)$.
We want to show that these vectors are linearly independent, so we have to show
\begin{equation}
\label{lindepend}
\sum_{\bs \in S(\lambda + \omega_i)} a_{\bs} f^{\bs}(v_{\lambda} \otimes v_{\omega_i}) = 0
\Rightarrow a_{\bs} = 0 \; \forall \; \bs \in S(\lambda + \omega_i).
\end{equation}
We may assume without loss of generality that all $\bs$ have the same weight,
say $\bs\in S(\lambda + \omega_i)^\mu$.
By Proposition~\ref{minimalset} we can split an element in $S(\lam+\omega_i)$ such that
$\bs = \bs_2 + \bm^{\bs}_i$, where $\bs_2 \in S(\lambda)$. Assume that we have a
non-trivial linear dependence relation in \eqref{lindepend}.
Fix $\bbs \in S(\lambda+ \omega_i)^\mu$ such that $a_{\bbs}\not=0$ in this relation
and $a_{\bt} = 0$ for all $\bt$ such that $\bm^{\bbs}_i \ll \bm^{\bt}_i$.
Consider first $\bbs = \bbs_2 + \bm^{\bbs}_i$, so we have
\begin{equation}
\label{rest}
f^{\bbs}(v_{\lambda} \otimes v_{\omega_i}) = c_{\bm^{\bbs}_i} f^{\bbs_2}v_{\lambda} \otimes f^{\bm^{\bbs}_i}v_{\omega_i} + \text{\it  other terms, }
\end{equation}
where $c_{\bm^{\bbs}_i}$ is a nonzero constant (product of binomial coefficients).

All the terms occurring in the linear dependence relation \eqref{lindepend} can be rewritten
as sums of terms of the form $f^{\br_2}v_{\lambda} \otimes f^{\br_1}v_{\omega_i}$. So in order to
prove that necessarily $a_\bs=0$ for all terms in \eqref{lindepend}, it is sufficient to show
that that the terms $f^{\br_2}v_{\lambda} \otimes f^{\br_1}v_{\omega_i}$ satisfying
$wt(\br_2)=wt(\bbs_2)$ and $wt(\br_1)=wt(\bm^{\bbs}_i)$ are linearly independent.

Let us first consider the possible terms in \eqref{rest} occurring among the \textit{other terms}. It is a sum of elements
$f^{\br_2}v_{\lambda} \otimes f^{\br_1}v_{\omega_i}$, where
$\br_2 + \br_1 = \bbs$ and $\br_1 \neq \bm^{\bbs}_i$. If $\wt(\br_1) = \wt(\bm^{\bbs}_i)$,
then either $\br_1\in S(\omega_i)$, but then $\br_1=\bm^{\bbs}_i$ for weight reasons, or
$\br_1\not\in S(\om_i)$.
In the latter case the entries in $\br_1$ are zero for all $\al_{k,\ell}\not\in R_i$ because
of the special choice of the ordering, and hence $\br_1$ has to be a mutation of $\bm^{\bbs}_i$.
Then by Proposition~\ref{mutation}, $\bm^{\bbs}_i \ll \bm^{\br_1 + \br_2}_i = \bm^{\bbs}_i$ which is a contradiction.
So the \textit{other terms} consist only of tensors of the form $f^{\br_2}v_{\lambda} \otimes f^{\br_1}v_{\omega_i}$,
where $\wt({\br_2}) \neq \wt({\bbs_2})$ and $\wt({\br_1})
\neq \wt({\bm^{\bbs}_i})$, hence for proving linear independence we can
neglect these terms.

To obtain a non-trivial linear combination such that $a_\bt  \neq 0$ for some $\bt\not=\bbs$,
one needs an element $\bt \in S(\lambda + \omega_i)^\mu$ which can be splitted
$\bt = \bt_2 + \bt_1$ such that $\wt({\bt_2}) = \wt({\bbs_2})$, $\wt({\bt_1}) =
\wt({\bm^{\bbs}_i})$, and $f^{\bt_2} v_{\lambda} \neq  0, f^{\bt_1} v_{\omega_i} \neq 0$.

Suppose that one has such a $\bt = \bt_2 +\bt_1 $ and $\bt_1 \notin S(\omega_i)$. By the same arguments as above,
$\bt_1$ is a mutation of $\bm^{\bbs}_i$ and hence by Proposition~\ref{mutation}, $\bm^{\bbs}_i \ll \bm^{\bt}_i$.
But in this case we have by assumption $a_{\bt} = 0$, contradicting the fact $a_\bt  \neq 0$.

It follows $\bt_1 \in S(\omega_i)$ and hence, by weight arguments, $\bt_1=\bm^{\bbs}_i$
and $\bt = \bt_2 + \bm^{\bs}_i$, where $\bt_2 \neq \bs_2$.

So if a term of the form $f^{\bt_2}v_{\lambda} \otimes f^{\bt_1}v_{\omega_i}$
$\wt({\bt_2}) = \wt({\bbs_2})$, $\wt({\bt_1}) =
\wt({\bm^{\bbs}_i})$ occurs in the linear
dependence relation \eqref{lindepend}, then necessarily
$\bt_1=\bm^{\bbs}_i$. Hence, by Proposition~\ref{minimalset}, $\bt_2 \in S(\lambda)$.
Since the possible $\bt_2$ are different from $\bbs$ and by induction the terms
$\{f^{\bt_2}v_{\lambda} \otimes f^{\bm^{\bbs}_i}v_{\omega_i}\mid\bt_2\in S(\lam)\}$ are linearly independent,
it follows $a_{\bbs} = 0$, contradicting the assumption $a_{\bbs} \neq 0$.

Summarizing, we have shown that for the order fixed at the beginning of the proof 
the $f^\bs v_{\lam+\om_i}$, $\bs\in S(\lam+\om_i)$,
are linearly independent and form a basis. This implies in particular that 
$\sharp S(\lam+\om_i)=\dim V(\lam+\om_i)$.
Now by Theorem~\ref{linearindipencetheorem} we know that the $f^\bs v_{\lam+\om_i}$, $\bs\in S(\lam+\om_i)$,
span $V(\lam+\om_i)$ for any chosen total order. So, for dimension reason, they also have to be linearly independent
for any chosen order.
\qed
\section{Proof of Theorem A and applications}
In this section we collect some immediate consequences of the constructions
in Sections 2 and 3. The proof of Theorem~\ref{numberofelements} shows:
\begin{cor}
\par\noindent
$$
\dim V(\la) =\# S(\la)=\,\text{number of integral points in the polytope}\, P(\lam).
$$
\end{cor}
By the defining inequalities (see \ref{polytopeequation}) for the polytope $P(\lam)$
it is obvious that for two dominant integral weights $\lam,\mu$ we have $P(\lam)+P(\mu)\subseteq P(\lam+\mu)$,
and hence for the integral points we have $S(\lam)+S(\mu)\subseteq S(\lam+\mu)$, too.
In fact, the reverse implication is also true:
\begin{prop}
$S(\la)+S(\mu)=S(\lam+\mu)$.
\end{prop}
\proof
Set $\nu=\lam+\mu$ and write $\nu=\sum k_i\om_i$ as a sum of fundamental weights.
Proposition~\ref{minimalset} provides an inductive procedure to write
an element $\bs$ in $S(\nu)$ as a sum $\bs=\sum_{i=1}^n\sum_{j=1}^{k_i} \bm_{i,j}$
such that $\bm_{i,j}\in S(\omega_i)$ for all $1\le i\le n$, $1\le j\le k_i$. This sum can be
reordered in such a way that $\bs=\bs^1+\bs^2$, $\bs^1\in S(\lam)$, $\bs^2\in S(\mu)$,
so $\bs\in S(\lam)+S(\mu)$.
\qed

\bigskip
As an interesting application we obtain a combinatorial
character formula for the representation $V(\la)$.
Let $P$ be the weight lattice and
for $\bs \in S(\lambda)$ define the weight
$$
\wt(\bs) := \sum_{1 \leq j \leq k \leq n} s_{j,k} \alpha_{j,k}.
$$
Let $S(\lam)^\mu$ be the subset of elements such that $\mu=\lam-\wt(\bs)$ and
let ${\tt S}(\lam)^\mu:=\#\{ \bs \in S(\la) \; | \; \mu=\la-\wt(\bs)\}$ be the number of elements of this set.
We obtain as a consequence of Theorem~\ref{basis}:
\begin{prop}
$$
char V(\la) =\sum_{\mu\in P}{\tt S}(\lam)_\mu e^\mu.
$$
\end{prop}
\bigskip
The big advantage of our approach is that it provides also a combinatorial formula
for the graded character. Recall that $gr V(\la)$ carries an additional grading
on each weight space $V(\la)^\mu$ of $V(\la)$:
$$
gr V(\la)^\mu=\bigoplus_{s\ge 0} gr_s V(\la)^\mu=\bigoplus_{s\ge 0} V(\la)^\mu_s/V(\la)^\mu_{s-1}.
$$
The graded character of the weight space is the polynomial
$$
p_{\lam,\mu}(q):=\sum_{s\ge 0}(\dim V(\la)^\mu_s/V(\la)^\mu_{s-1})q^s
$$
and the graded character of $V(\lam)$ is
$$
char_q (V(\la))= \sum_{\mu\in P}p_{\lam,\mu}(q) e^\mu.
$$
We have a natural notion of a {\it degree} for the
multi-exponents:
\begin{dfn}
$$\deg(\bs) := \sum_{1 \leq j \leq k \leq n} s_{j,k}.$$
\end{dfn}
As an immediate consequence of Theorem~\ref{basis} we get
\begin{cor*}
$p_{\lam,\mu}(q)=\sum_{\bs\in S(\lam)^\mu} q^{\deg \bs}$ and
$$
char_q(V(\la)) = \sum_{\bs \in S(\la)} e^{\la-\wt(\bs)} q^{\deg(\bs)}.
$$
\end{cor*}

Finally, we note that the results of Sections \ref{span} and \ref{li}
imply the description of the annihilating ideal $I(\la)$.
\begin{thm}\label{idealtheorem}
\begin{equation}\label{ideal}
I(\la)=S(\n^-)\left(\U(\n^+)\circ \mathrm{span}\{f_\al^{(\la,\al)+1}, \al>0\}\right).
\end{equation}
\end{thm}
\begin{proof}
Since $f_\al^{(\la,\al)+1}v_\la=0$ in $V(\la)$ for all positive roots $\al$, the right hand side
of \eqref{ideal} belongs to $I(\la)$. Section \ref{span}
shows that the relations in the RHS of \eqref{ideal} are enough to rewrite
any element of $gr V(\la)$ in terms of the basis element $f^\bs v_\la$, $\bs\in S(\la)$.
This proves our theorem.
\end{proof}

\section*{Acknowledgements}
The work of Evgeny Feigin was partially supported
by the Russian President Grant MK-281.2009.1, the RFBR Grants 09-01-00058, 07-02-00799
and NSh-3472.2008.2, by Pierre Deligne fund
based on his 2004 Balzan prize in mathematics and by Alexander von
Humboldt Fellowship.
The work of Ghislain Fourier was partially supported by the DFG project
``Kombinatorische Beschreibung von Macdonald und Kostka-Foulkes Polynomen``.
The work of Peter Littelmann was partially supported by the
priority program SPP 1388 of the German Science Foundation.

\end{document}